\pdfoutput=1

\documentclass[11pt,a4paper]{article}

\usepackage{natbib}
\usepackage{anysize}
\usepackage{color}
\usepackage{graphicx}

\usepackage[mathscr]{euscript}

\usepackage{amsmath,amssymb,amsthm,amsopn}

\usepackage{pgfplots}

\pgfplotsset{compat=newest}

\theoremstyle{remark}

\newcommand{\mat}[1]{\mathbf{#1}}
\newcommand{\D}{\mathsf{D}}
\newcommand{\bH}{{\bf H}}

\newcommand{\bx}{{\boldsymbol x}}
\newcommand{\by}{{\boldsymbol y}}

\newcommand{\bX}{{\bf X}}

\newcommand{\bmu}{{\boldsymbol \mu}}

\newcommand{\bSigma}{{\boldsymbol \Sigma}}

\newcommand{\bgamma}{{\boldsymbol \gamma}}
\newcommand{\bnu}{{\boldsymbol \nu}}

\theoremstyle{plain}

\newtheorem{teor*}{Teorema}

\theoremstyle{definition}

\title{Clusters and water flows: a novel approach to modal clustering through Morse theory}

\author{Jos\'e E. Chac\'on\footnote{Departamento de
Matem\'aticas, Universidad de Extremadura, E-06006 Badajoz, Spain. E-mail:
{\tt jechacon@unex.es}}}

\begin{document}

\maketitle

\begin{abstract}
\noindent The problem of finding groups in data (cluster analysis) has been extensively studied by researchers from the fields of Statistics and Computer Science, among others. However, despite its popularity it is widely recognized that the investigation of some theoretical aspects of clustering has been relatively sparse. One of the main reasons for this lack of theoretical results is surely the fact that, unlike the situation with other statistical problems as regression or classification, for some of the cluster methodologies it is quite difficult to specify a population goal to which the data-based clustering algorithms should try to get close. This paper aims to provide some insight into the theoretical foundations of the usual nonparametric approach to clustering, which understands clusters as regions of high density, by presenting an explicit formulation for the ideal population clustering.

\end{abstract}

%\medskip
%\noindent {\em Keywords: cluster analysis, Morse Theory}

\newpage

\section{Introduction}\label{sec:1}

Clustering is one of the branches of Statistics with more research activity in the recent years. As
noted by \cite{M07}, ``clustering is a young domain of research, where rigorous methodology is
still striving to emerge". Indeed, some authors have recently expressed their concerns about the
lack of theoretical or formal developments for clustering \citep{vLBD05,ABD08,ZBD09}. This paper
aims to contribute to this regularization (or, say, rigorousization).

Stated in its most simple form, cluster analysis consists in ``partitioning a data set into groups so
that the points in one group are similar to each other and are as different as possible from the points in other groups" \citep[p. 293]{HMS01}. Posed as such, the problem does not even seem to have a statistical meaning.

In fact, some clustering techniques are solely based on the distances between the observations. Close observations are joined together to form a group, and extending the notion of inter-point distance to distance between groups, the resulting groups are gradually merged until all the initial observations are contained into a single group. This represents, of course, the notion of agglomerative {\it hierarchical clustering} \citep[Section 12.3]{I08}. The graphical outcome depicting the successive agglomeration of data points up to a single group is the well-known dendrogram, and depending on the notion of inter-group distance used along the merging process the most common procedures of this type are known as single linkage, complete linkage or average linkage \citep[see also][p. 523]{HTF09}.

A first statistical flavour is noticed when dealing with those clustering methodologies that represent each cluster by a central point, such as the mean, the median or, more generally, a trimmed mean. This class of techniques is usually referred to as {\it partitioning methods}, and surely the most popular of its representatives is $K$-means \citep{MQ67}. For a pre-specified number $K$ of groups, these algorithms seek for $K$ centers with the goal of optimizing a certain score function representing the quality of the clustering  \citep[Chapter 5]{Eal11}.

A fully statistical notion of the problem is {\it distribution-based clustering}. This approach
further assumes that the data set is a sample from a probability distribution. As with all the
statistical procedures, there are parametric and nonparametric methodologies for distribution-based
clustering. Surely the golden standard of parametric clustering is achieved through mixture
modeling as clearly described in \cite{FR02}. It is assumed that the distribution generating the
data is a mixture of simple parametric distributions, e.g. multivariate normal distributions, and
each component of the mixture is associated to a different population cluster. Maximum likelihood
is used to fit a mixture model and then each data point is assigned to the most likely component
using the Bayes rule.

Eventually, by using mixture distribution components defined through Mahalanobis distances, some partitioning methods and mixture-model clustering can be regarded to have the same motivation, and under this point of view some partitioning techniques have evolved to reach a high level of maturity, leading to very sophisticated methods. See, for instance, the recent review of \cite{GEal10} and references therein.

In parametric distribution-based clustering the notion of population cluster is clearly defined through the components of the mixture model. If the true underlying distribution belongs indeed to the postulated mixture model, this identification is certainly quite adequate. However, often mixture models are used in a ``nonparametric way", on the basis that they usually represent a  dense subclass of the set of all continuous distributions and, therefore, they are useful to approximate arbitrarily well any continuous distribution. In that case, the identification of components with clusters may not coincide with the natural intuition. For instance, several normal components may be needed to model accurately a distribution with a heavier-than-normal tail, yet the resulting normal mixture might as well be a unimodal density, suggesting the existence of a single group, though with a non-normal distribution shape. The fact that mixture components and clusters are not necessarily the same has motivated several variations of mixture-model clustering to amend this deficiency; see, for example, \cite{BCG00}, \cite{Bal10} or \cite{H10}.

The previous example is helpful to understand the motivation of the nonparametric methodology, which is based on identifying clusters as regions of high density separated from each other by regions of lower density. In this sense, population clusters are naturally associated with the modes (i.e., local maxima) of the
probability density function, and this nonparametric approach is denominated mode-based clustering or {\it modal clustering} \citep{LRL07}.  Precisely, each cluster is usually understood as the ``domain of attraction" of a mode \citep{S03}.

The concept of domain of attraction is not that simple to specify, and providing a precise definition is one of the goals of this paper that will be explored in the next section. It is perhaps this imprecise notion of population cluster in the nonparametric case what does not fully convince many researchers and statisticians and lead them to adopt one of the alternative previously expounded techniques (parametric clustering, partitioning methods or hierarchical algorithms).

The first attempt to make the goal of modal clustering precise is through the notion of level sets
\citep{H75}. If the distribution of the data has a density $f$, given $\lambda\geq0$ the
$\lambda$-level set of $f$ is defined as $L(\lambda)=\{\bx\colon f(\bx)\geq\lambda\}$. Then,
population $\lambda$-clusters are defined as the connected components of $L(\lambda)$, a definition
that clearly captures the notion of groups having a high density. An extensive account of the
usefulness of level sets in applications is given in \cite{MP09}.

The main drawback of this approach is the fact that the notion of population cluster depends on the
level $\lambda$, as recognized by \cite{S03}. However, other authors like \cite{CFF01} or
\cite{CPP11}, affirm that the choice of $\lambda$ is only a matter of resolution level of the
analysis. Indeed, the level sets can be re-parameterized in terms of their probability content, as
in \cite{C06} or \cite{AT07}, for ease of interpretability. This is usually done by finding, for
any $p\in(0,1)$, the threshold $\lambda_p$ defined as the largest value of $\lambda$ such that the
probability of $L(\lambda)$ is greater than $p$, and setting $R(p)=L(\lambda_p)$ so that the level
set $R(p)$ has probability $p$ as long as the border of $R(p)$ has null probability. Then, the
chosen value of $p$ implies that a proportion $\alpha=1-p$ of the data is left out of the
clustering mechanism, and therefore the level $\lambda$ (or $p$, or $\alpha$) can be interpreted as
a trimming parameter motivated by robustness considerations. In fact, the parametrization of the
level sets in terms of the probability $\alpha$ of the set of non-classified observations is also
commonly used; see \cite{H96,BCC01} or \cite{SW10}.

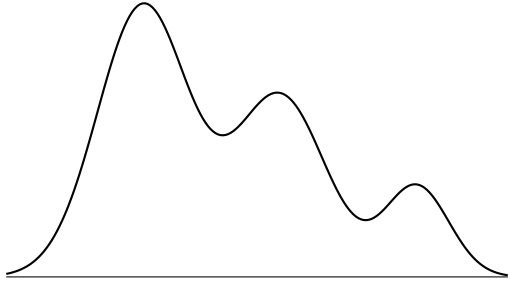
\begin{figure}[t]\centering
\begin{tikzpicture}[scale=0.6]
    \def\normalmixt{\x,{6*exp(-((\x)^2)/2)+4*exp(-((\x-3)^2)/2)+2*exp(-((\x-6)^2))}}
    \draw[color=black,domain=-3:8,samples=200,thick] plot (\normalmixt) node[right] {};
    \draw (-3,0) -- (8,0) node[right] {};
\end{tikzpicture}
\caption{Univariate trimodal density for which it is not possible to capture its whole cluster structure using a level set analysis based on a single level.}
\label{fig:1}
\end{figure}

Nevertheless, it is easy to think of many examples in which it is impossible to observe the whole cluster structure on the basis of a single level $\lambda$.
Essentially as in \citet[][p. 906]{RSNW12}, Figure \ref{fig:1} shows a simple univariate example of this phenomenon: three different groups are visually identifiable, yet none of the level sets of the density has three connected components. To amend this, the usual recommendation is to analyze the cluster structure for several values of the level $\lambda$ (or the coverage probability $p$). Graphical tools oriented to this goal are the cluster tree \citep{S03}, or the mode function \citep{AT07}. Both graphics are useful to show how the clusters emerge as a function of $\lambda$ or $p$, respectively. A detailed description of the cluster identification methodology through the cluster tree is given in Section \ref{sec:2} below.

Level sets also arise as the maximizers of excess mass functionals. For a given value of
$\lambda\geq0$, the excess mass functional is defined as $E_{\lambda}(A)=P(A)-\lambda\mu(A)$ for
any measurable set $A$, where $\mu$ denotes the Lebesgue measure. \cite{H87} and \cite{MS87,MS91}
independently introduced this functional and noticed that it reaches its maximum value precisely at
$L(\lambda)$. A nice feature of the excess mass functional is that it admits a fairly simple
empirical estimator obtained by replacing $P$ for the empirical probability measure of the sample.
However, the maximizer of this empirical criterion does not provide a useful estimate of
$L(\lambda)$ because it turns out to be just the set of all the data points. Nevertheless, by
restricting the empirical functional to a suitable class of measurable sets its maximizer can be
shown to be a sensible estimate of $L(\lambda)$ within such a class. This estimation procedure is
studied in detail in \cite{P95} making use of empirical processes theory.

On the other hand, the estimation of density level sets using smoothing techniques has its roots in
the papers of \cite{CF97,T97} or \cite{W97}, and related data-driven cluster methodologies based on
the notion of $\lambda$-clusters are proposed in \cite{CFF01,S03,AT07,NS10,SN10,RW10} or
\cite{RSNW12}, among others. In all these papers a plug-in approach is adopted, in the sense that a
population level set is estimated by the corresponding level set of some density estimator. It is
to be noted that for most of the plug-in estimators, and also for the excess mass estimators, the
computation of empirical level sets is indeed a difficult issue, and many of the aforementioned
references also present intricate algorithms to approximate the corresponding level set estimators.
Relatedly, in Machine Learning literature surely the most popular algorithms for density-based
clustering include DBSCAN \citep{Eal96}, based on a single level, and OPTICS \citep{Aal99}, which
is a generalization of the former which allows to compute the cluster tree induced by a range of
levels.

\section{A precise approach to the goal of modal clustering}\label{sec:2}

Certainly, the very concept of modal cluster may admit some well-supported criticism, as
illustrated in the thorough reflections of \citet[][Section 6]{S03} or \cite{H10}, and we
agree with \citet[][p. 408]{I08} that ``there is no universally accepted definition of exactly what
constitutes a cluster". Some of the most common definitions have been surveyed in the previous
section. However, it is not the goal of this paper to contribute to the discussion of the different
notions of cluster, but to introduce a precise definition of population modal cluster which does
not rely on the concept of level set. This aims to provide an answer to Question 1 in \cite{vLBD05}: ``How does a desirable clustering look like if we have complete knowledge about our data generating process?"

Since the empirical formulation of the clustering task consists of partitioning a data set into
homogeneous groups, a population clustering $\mathscr C$ of a probability distribution $P$ on $\mathbb R^d$ should
be understood as a partition of $\mathbb R^d$ into mutually disjoint measurable components
$C_1,\dots,C_k$, each with positive probability content \citep{BLP06}. %,CM10}.
The components of such a partition are called population clusters. Thus, two population clusterings
$\mathscr C$ and $\mathscr D$ are identified to be the same if they have the same number of
clusters and, up to a permutation of the cluster labels, every cluster of $\mathscr C$ and its
most similar match in $\mathscr D$ differ in a null-probability set (more details on this are elaborated
in Section \ref{sec:cons}).

As anticipated in the previous section, a basic principle about the concept of ideal population clustering is that it should reflect the notion of a
partition into regions of high density separated from each other by regions of lower density, leading to the idea of population cluster as the domain of
attraction of a density mode. We begin our way to formalize this concept by
showing some illustrative examples in one and two dimensions.

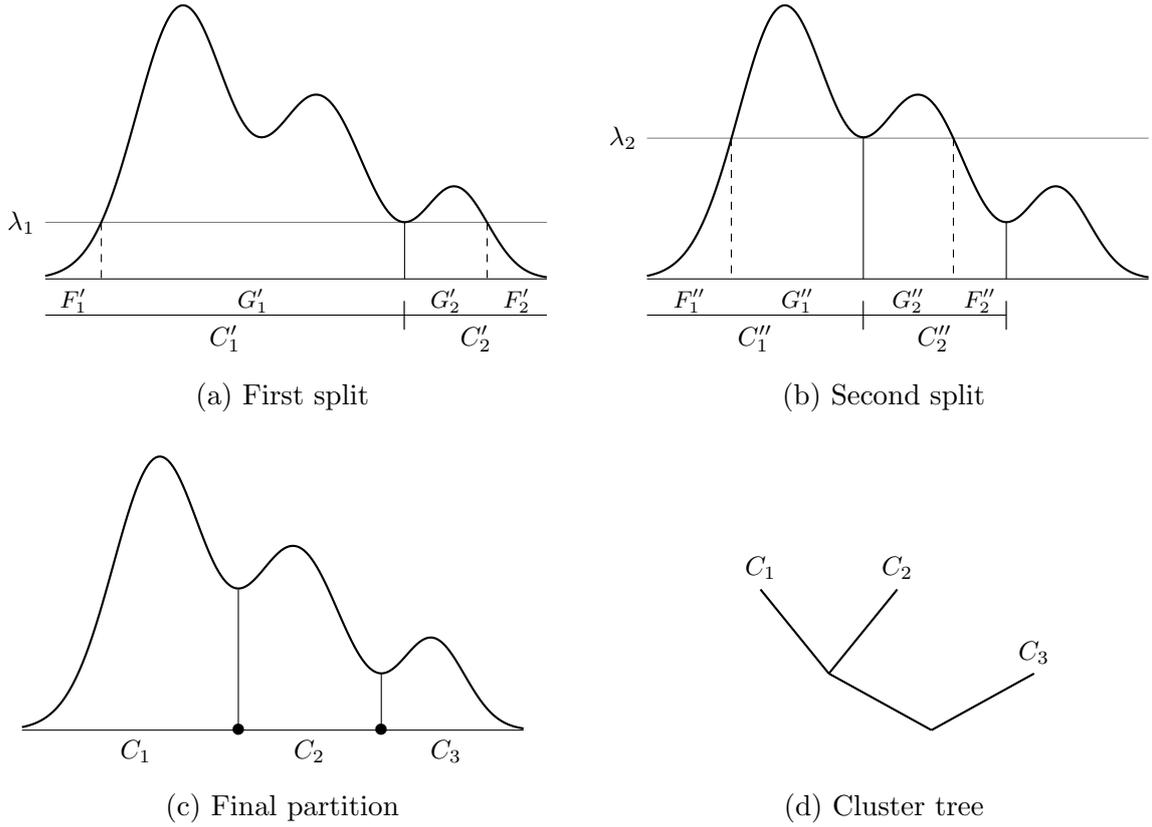
\begin{figure}[t]\centering

\begin{tabular}{@{}cc@{}}

\begin{tikzpicture}[scale=0.6]
    \def\normalmixt{\x,{6*exp(-((\x)^2)/2)+4*exp(-((\x-3)^2)/2)+2*exp(-((\x-6)^2))}};
    \draw[color=gray,thin] (-3,1.25) -- (8,1.25) node[right] {};
    \draw[color=black,domain=-3:8,samples=200,thick] plot (\normalmixt) node[right] {};
    \draw (-3,0) -- (8,0) node[right] {};
    \draw (-3,1.25) node[left] {\small $\lambda_1$};
    \draw[] (4.88,1.25) -- (4.88,0) node[below] {};
    \draw[dashed] (-1.77,1.25) -- (-1.77,0) node[below] {};
    \draw[dashed] (6.69,1.25) -- (6.69,0) node[below] {};

    \draw[] (-3,-.8) -- (8,-.8) node[] {};
    \draw (4.88,-0.8) node[] {$|$};
    \draw(-2.38,0) node[below] {\footnotesize $F_1'$};
    \draw(1.55,0) node[below] {\footnotesize $G_1'$};
    \draw(5.78,0) node[below] {\footnotesize $G_2'$};
    \draw(7.34,0) node[below] {\footnotesize $F_2'$};

    \draw (0.94,-.8) node[below] {\small $C_1'$};
    \draw (6.44,-.8) node[below] {\small $C_2'$};

\end{tikzpicture}&\begin{tikzpicture}[scale=0.6]
    \def\normalmixt{\x,{6*exp(-((\x)^2)/2)+4*exp(-((\x-3)^2)/2)+2*exp(-((\x-6)^2))}};
    \draw[color=gray,thin] (-3,3.11) -- (8,3.11) node[right] {};
    \draw[color=black,domain=-3:8,samples=200,thick] plot (\normalmixt) node[right] {};
    \draw (-3,0) -- (8,0) node[right] {};
    \draw (-3,3.11) node[left] {\small $\lambda_2$};
    \draw[] (1.745,3.11) -- (1.745,0) node[below] {};
    \draw[] (4.88,1.25) -- (4.88,0) node[below] {};
    \draw[dashed] (-1.15,3.11) -- (-1.15,0) node[below] {};
    \draw[dashed] (3.72,3.11) -- (3.72,0) node[below] {};

    \draw[] (-3,-.8) -- (4.88,-.8) node[] {};
    \draw (4.88,-0.8) node[] {$|$};
    \draw (1.745,-0.8) node[] {$|$};
    \draw(-2.07,0) node[below] {\footnotesize $F_1''$};
    \draw(0.3,0) node[below] {\footnotesize $G_1''$};
    \draw(2.73,0) node[below] {\footnotesize $G_2''$};
    \draw(4.3,0) node[below] {\footnotesize $F_2''$};

    \draw (-0.63,-.8) node[below] {\small $C_1''$};
    \draw (3.31,-.8) node[below] {\small $C_2''$};

\end{tikzpicture}\\

(a) First split & (b) Second split\\[.5cm]

\begin{tikzpicture}[scale=0.6]
    \def\normalmixt{\x,{6*exp(-((\x)^2)/2)+4*exp(-((\x-3)^2)/2)+2*exp(-((\x-6)^2))}};
    \draw[color=black,domain=-3:8,samples=200,thick] plot (\normalmixt) node[right] {};
    \draw (-3,0) -- (8,0) node[right] {};

    \draw[] (1.745,3.11) -- (1.745,0) node[below] {};
    \draw[] (4.88,1.25) -- (4.88,0) node[below] {};

    \draw (4.88,0) node[] {$\bullet$};
    \draw (1.745,0) node[] {$\bullet$};
    \draw (-0.5,0) node[below] {\small $C_1$};
    \draw (3.3,0) node[below] {\small $C_2$};
    \draw (6.3,0) node[below] {\small $C_3$};
\end{tikzpicture}&\begin{tikzpicture}[scale=0.6]
    \draw[color=white] (-3,0) -- (8,0) node[right] {};
    \draw[thick] (3.75,0) -- (6,1.25) node[] {};
    \draw[thick] (3.75,0) -- (1.5,1.25) node[] {};

    \draw[thick] (1.5,1.25) -- (0,3.11) node[] {};
    \draw[thick] (1.5,1.25) -- (3,3.11) node[] {};

    \draw (0,3.11) node[above] {\small $C_1$};
    \draw (6,1.25) node[above] {\small $C_3$};
    \draw (3,3.11) node[above] {\small $C_2$};

    \draw[color=white] (-0.5,0) node[below] {\small $C_1$};

\end{tikzpicture}\\

(c) Final partition& (d) Cluster tree\\

\end{tabular}

\caption{Identification of clusters for the trimodal density example using the cluster tree.}
\label{fig:2}
\end{figure}

In the one-dimensional case, it seems clear from Figure \ref{fig:2} how this can be achieved. To
begin with, the level set methodology identifies the three clusters in the density depicted in
Figure \ref{fig:1} by computing the cluster tree as described clearly in \cite{NS10}: starting from
the $0$-level set, which corresponds to the whole real line in this example and consists of a single connected component, $\lambda$ is increased until it reaches $\lambda_1$, where two
components for the $\lambda_1$-level set are found, $G_1'$ and $G_2'$, resulting in the cluster
tree splitting into two different branches (see Fig. \ref{fig:2}, panel (a)). These two components
$G_1'$ and $G_2'$ are usually called cluster cores. They do not constitute a partition of $\mathbb R$, but the
remaining parts $F_1'$ and $F_2'$, referred to as fluff in \cite{NS10}, can be assigned to either
the left or the right branch depending on whichever of them is closer. Thus, at level $\lambda_1$
the partition $\mathbb R=C_1'\cup C_2'$ is obtained. The point dividing the line into these two
components can be arbitrarily assigned to either of them; this assignment makes no difference,
because it leads to equivalent clusterings.

At level $\lambda_2$ the left branch is further divided into two branches (see panel (b) of Fig. \ref{fig:2}). Again, the two cluster core components $G_1''$ and $G_2''$ do not form a partition of the set $C_1'$ associated with the
previous node of the tree, but it is clear how the fluff $F_1''$ and $F_2''$ can be assigned to form a
partition $C_1''\cup C_2''$ of $C_1'$. Since no further splitting of the cluster tree is observed as
$\lambda$ increases, the final partition of $\mathbb R$ is $C_1''\cup C_2''\cup C_2'$, renamed to $C_1\cup
C_2\cup C_3$ in panel (c) of Figure \ref{fig:2}.

It is immediate to observe that the levels at which a connected component breaks into two different
ones correspond precisely to local minima of the density function, so an equivalent formulation
consists of defining population clusters as the connected components of $\mathbb R$ minus the points where a local minimum is attained (the solid circles in panel (c) of Figure \ref{fig:2}). Notice that this definition does not involve the computation of level sets for a range of levels, nor their cores and fluff, so that it constitutes a more straightforward approach to the same concept in the unidimensional setup.

\pgfplotsset{ colormap={whitered}{color(0cm)=(white); color(1cm)=(cyan!60!blue)} }

\begin{figure}[t]\centering

\begin{tabular}{@{}cc@{}}

\begin{tikzpicture}[
    scale=0.8,
    declare function={mu1=1;},
    declare function={mu2=2;},
    declare function={sigma1=0.5;},
    declare function={sigma2=1;},
    declare function={normal(\m,\s)=1/(2*\s*sqrt(pi))*exp(-(x-\m)^2/(2*\s^2));},
    declare function={norm(\ma,\sa,\mb,\sb)=
        1/(2*pi*\sa*\sa) * exp(-((-\ma)^2/\sa^2 + (x-\mb)^2/\sb^2))/2;},
    declare function={bivar(\ma,\sa,\mb,\sb)=
        1/(2*pi*\sa*\sa) * exp(-((x-\ma)^2/\sa^2 + (y-\mb)^2/\sb^2))/2;}]
\begin{axis}[
    colormap name=whitered,
    width=.6\textwidth,
    view={35}{45},
    enlargelimits=false,
    grid=major,
    domain=-3:3,
    y domain=-3:3,
    ytick={-3,-2,-1,0,1,2,3},
    ztick={0,0.02,0.04,0.06},
    samples=40%%40
]

\addplot3 [surf] {bivar(-1.5,1,0,1)+bivar(1.5,1,0,1)};

\addplot3 [domain=-3:0.3,samples=200, samples y=0, very thick, smooth]
(0,x,{norm(-1.5,1,0,1)+norm(1.5,1,0,1)});

\addplot3 [domain=0.3:3,samples=200, samples y=0, smooth, thick, color=gray]
(0,x,{norm(-1.5,1,0,1)+norm(1.5,1,0,1)});

\node at (axis cs:0,0,0.01677481) {$\bullet$};
\node[above] at (axis cs:0,0,0.01677481) {$\bx_0$};

\node at (axis cs:-1.5,0,0.07958729) {$\bullet$};
\node[above] at (axis cs:-1.5,0,0.07958729) {$\bx_1$};

\node at (axis cs:1.5,0,0.07958729) {$\bullet$};
\node[above] at (axis cs:1.5,0,0.07958729) {$\bx_2$};

\end{axis}
\end{tikzpicture}&\includegraphics[width=0.4\textwidth]{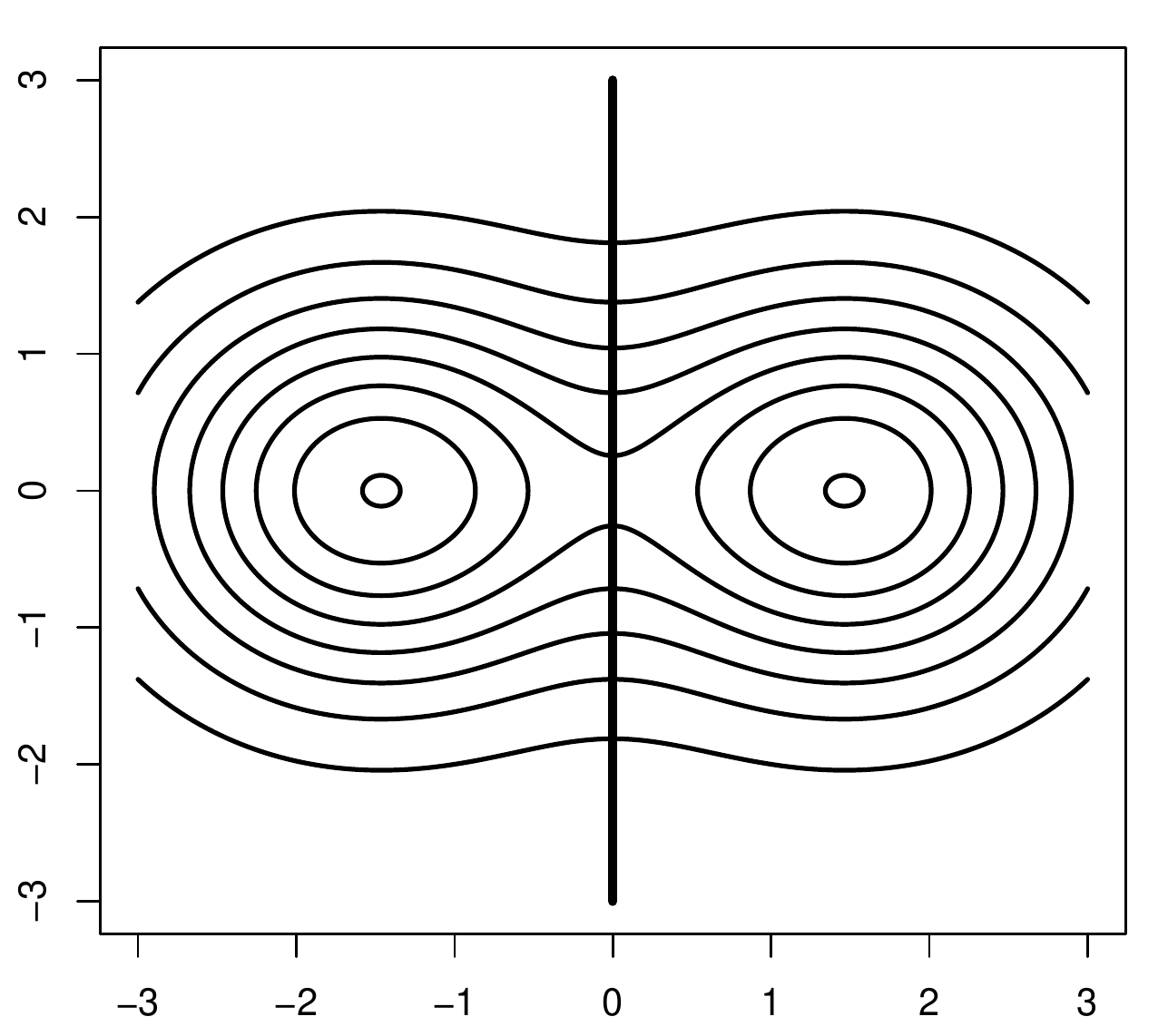}
%\begin{tikzpicture}[
%    scale=0.8,
%        declare function={bivar(\ma,\sa,\mb,\sb)=
%        1/(2*pi*\sa*\sa) * exp(-((x-\ma)^2/\sa^2 + (y-\mb)^2/\sb^2))/2;}]
%    \begin{axis}[width=0.6\textwidth,view={0}{90},axis equal, colormap/blackwhite,domain=-3:3,xtick={-3,-2,-1,0,1,2,3},
%        enlargelimits=false,ytick={-3,-2,-1,0,1,2,3}]%,enlarge x limits]
%    \addplot3[id=cont,contour gnuplot={labels=false,number=12},thick,samples=40]
%        {bivar(-1.5,1,0,1)+bivar(1.5,1,0,1)};
%    \addplot3[color=black,ultra thick] coordinates {(0,-3,0) (0,3,0)};
%
%\node at (axis cs:0,0,0.01677481) {$\bullet$};
%\node[above right] at (axis cs:0,0,0.01677481) {$\!\bx_0$};
%
%\node at (axis cs:-1.5,0,0.07958729) {$\bullet$};
%\node[above] at (axis cs:-1.5,0,0.07958729) {$\bx_1$};
%
%\node at (axis cs:1.5,0,0.07958729) {$\bullet$};
%\node[above] at (axis cs:1.5,0,0.07958729) {$\bx_2$};
%
%    \end{axis}
%
%    \draw[color=white,ultra thick] (-.1,-1) -- (.1,-1);
%
%\end{tikzpicture}
\end{tabular}
\caption{Bidimensional example, with two groups clearly identifiable at an intuitive level.}
\label{fig:3}
\end{figure}
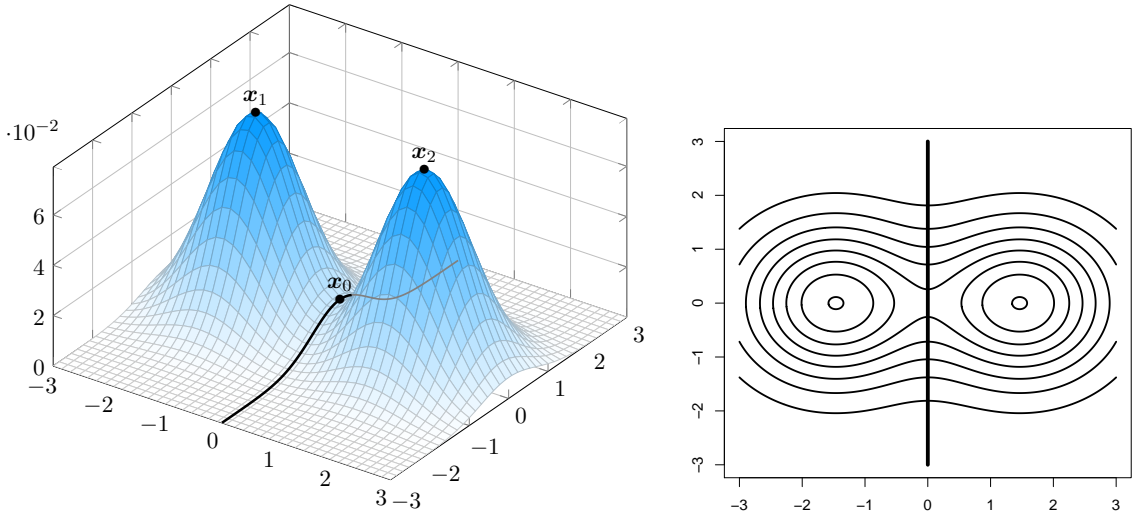

To get an idea of how to generalize the previous approach to higher dimensions, consider the
following extremely simple bidimensional example: an equal-proportion mixture of two normal
distributions, each with identity variance matrix and centered at $\bmu_1=(-\frac32,0)$ and
$\bmu_2=-\bmu_1$, respectively. At an intuitive level, it is clear from Figure \ref{fig:3} that the
most natural border to separate the two visible groups is the black line. The problem is: what is
exactly that line? Is it identifiable in terms of the features of the density function in a
precise, unequivocal way? A nice way to answer to these questions is by means of Morse theory.

Morse theory is a branch of Differential Geometry that provides tools for analyzing the topology of
a manifold $M\subseteq\mathbb R^d$ by studying the critical points of a differentiable function
$f\colon M\to\mathbb R$. A classical reference book on this subject is \cite{Mil63} and enjoyable
introductions to the topic can be found in \cite{Mat02} and \citet[][Chapter 7]{Jos11}. A useful
application of Morse theory is for terrain analysis, as nicely developed in \cite{Vit10}. In
terrain analysis, a mountain range can be regarded as the graph of a function $f\colon M\to\mathbb
R$, representing the terrain elevation, over a terrain $M\subseteq\mathbb R^2$, just as in the left
graphic of Figure \ref{fig:3}. The goal of terrain analysis is to provide a partition of $M$
through watersheds indicating the different regions, or catchment basins, where water flows under
the effect of gravity.

The fundamentals of Morse theory can be extremely summarized as follows. A smooth function $f\colon
M\to\mathbb R$ is called a {\it Morse function} if all its critical points are non-degenerate. Here
the critical points of $f$ are understood as those $\bx_0\in M$ for which the gradient $\D
f(\bx_0)$ is null, and non-degeneracy means that the determinant of the Hessian matrix $\mathsf{H}
f(\bx_0)$ is not zero. For such points the {\it Morse index} $m(\bx_0)$ is defined as the number of
negative eigenvalues of $\mathsf{H}f(\bx_0)$.

Morse functions can be expressed in a fairly simple form in a neighborhood of a critical point $\bx_0$, as the result known as Morse lemma shows that it is possible to find local coordinates $x_1,\dots,x_n$ such that $f$ can be written as $\pm x_1\pm\cdots\pm x_d+c$ around $\bx_0$, where the number of minus signs in the previous expression is precisely $m(\bx_0)$. For example, for $d=2$ the three possible configurations for a critical point are shown in Figure \ref{fig:4}, corresponding to a local minimum, a saddle point and a local maximum (from left to right), with Morse indexes 0, 1 and 2, respectively.

\pgfplotsset{ colormap={whitered}{color(0cm)=(white); color(1cm)=(cyan!60!blue)} }
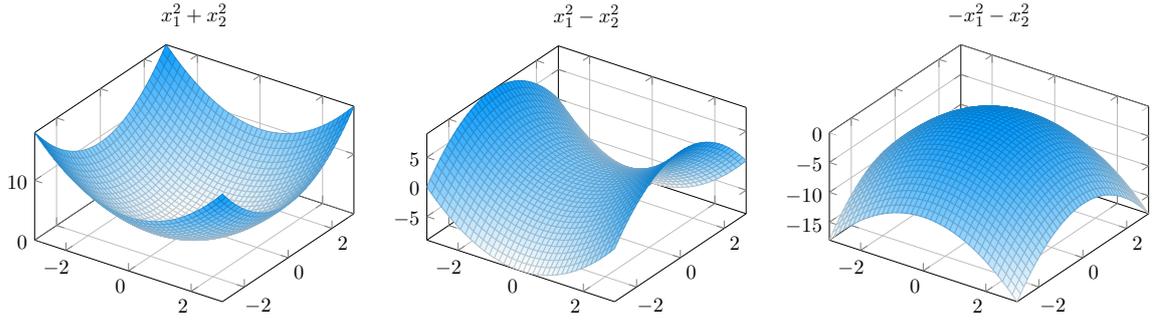
\begin{figure}[t]\centering

\begin{tabular}{@{}ccc@{}}

\begin{tikzpicture}[
    scale=0.7,
    declare function={pmax=x^2+y^2;},
    declare function={pmin=-x^2-y^2;},
    declare function={psilla=x^2-y^2;}
    ]
\begin{axis}[
    colormap name=whitered,
    width=0.38\textwidth,
    view={35}{45},
    enlargelimits=false,
    grid=major,
    domain=-3:3,
    y domain=-3:3,
    samples=40,%0,
    %xlabel=$x_1$,
    %ylabel=$x_2$,
    title={$x_1^2+x_2^2$}
]

\addplot3 [surf] {pmax};
\end{axis}
\end{tikzpicture}&

\begin{tikzpicture}[
    scale=0.7,
    declare function={pmax=x^2+y^2;},
    declare function={pmin=-x^2-y^2;},
    declare function={psilla=x^2-y^2;}
    ]
\begin{axis}[
    colormap name=whitered,
    width=0.38\textwidth,
    view={35}{45},
    enlargelimits=false,
    grid=major,
    domain=-3:3,
    y domain=-3:3,
    samples=40,%0,
    %xlabel=$x_1$,
    %ylabel=$x_2$,
    title={$x_1^2-x_2^2$}
]

\addplot3 [surf] {psilla};
\end{axis}
\end{tikzpicture}&

\begin{tikzpicture}[
    scale=0.7,
    declare function={pmax=x^2+y^2;},
    declare function={pmin=-x^2-y^2;},
    declare function={psilla=x^2-y^2;}
    ]
\begin{axis}[
    colormap name=whitered,
    width=0.38\textwidth,
    view={35}{45},
    enlargelimits=false,
    grid=major,
    domain=-3:3,
    y domain=-3:3,
    samples=40,%0,
    %xlabel=$x_1$,
    %ylabel=$x_2$,
    title={$-x_1^2-x_2^2$}
]

\addplot3 [surf] {pmin};
\end{axis}
\end{tikzpicture}

\end{tabular}
\caption{The three possible configurations around a critical point of a Morse function in the bidimensional case.}
\label{fig:4}
\end{figure}

The decomposition of $M$ suggested by Morse theory is made in terms of the unstable and/or stable manifolds of the critical points of $f$ as explained next. Consider the initial value problem defined by the minus gradient vector of $f$. For a given value of $\bx\in M$ at time $t=0$, the integral curve $\bnu_\bx\colon\mathbb R\to M$ of such initial value problem is the one satisfying
\begin{equation*}
\bnu_\bx'(t)=-\D f\big(\bnu_\bx(t)\big),\  \bnu_\bx(0)=\bx
\end{equation*}
and the set of all these integral curves is usually referred to as the negative gradient flow. Since the minus gradient vector defines the direction of steepest descent of $f$, these curves (or properly speaking, their images through $f$) represent the trajectories of the water flow subject to gravity.

With respect to the negative gradient flow, the {\it unstable manifold} of a critical point $\bx_0$ is defined as the set of points whose integral curve starts at $\bx_0$, that is
\begin{equation*}
W^u_-(\bx_0)=\big\{\bx\in M\colon \lim_{t\to-\infty}\bnu_\bx(t)=\bx_0\big\}.
\end{equation*}
Analogously, the stable manifold of $\bx_0$ is the set of points whose integral curve finishes at $\bx_0$, that is, $W^s_-(\bx_0)=\big\{\bx\in M\colon \lim_{t\to+\infty}\bnu_\bx(t)=\bx_0\big\}$. It was firstly noted by \cite{Thom49} that the set of unstable manifolds corresponding to all the critical points of $f$ provides a partition of $M$ (the same is true for the stable manifolds). Furthermore, the unstable manifold $W^u_-(\bx_0)$ has dimension $m(\bx_0)$.

The main contribution of this section is to define the population clusters of a density $f$ as the unstable
manifolds of the negative gradient flow corresponding to local maxima of $f$. Or in a more prosaic
way, in terms of water flows, a cluster is just the region of the terrain that would be flooded by
a fountain emanating from a peak of the mountain range.

Although this is an admittedly cumbersome definition, going back to Figure \ref{fig:3} it is clear
that it just describes the notion that we were looking for. The critical point $\bx_0=(0,0)$ is a
saddle point, thus having Morse index 1, and the black line is precisely its associated unstable
manifold, $W^u_-(\bx_0)=\{0\}\times\mathbb R$, which is a manifold of dimension 1. The remaining
two critical points are local maxima, and their respective unstable manifolds are
$W^u_-(\bx_1)=(-\infty,0)\times\mathbb R$ and $W^u_-(\bx_2)=(0,\infty)\times\mathbb R$, manifolds
of dimension 2 so that we can partition $\mathbb R^2=W^u_-(\bx_0)\cup W^u_-(\bx_1)\cup
W^u_-(\bx_2)$, showing $W^u_-(\bx_1)$ and $W^u_-(\bx_2)$ as two population clusters separated by
the border $W^u_-(\bx_0)$, which is a null-probability set.

Notice that this definition also applies to the previous univariate example in Figure \ref{fig:2}: the clusters $C_1$, $C_2$ and $C_3$ are just the unstable manifolds of the three local maxima (they are manifolds of dimension 1), and for the two local minima their unstable manifolds have dimension 0, so they include only the respective points of minima.

Moreover, if we focus on the gradient flow, instead of the negative gradient flow, then its integral curves satisfy
\begin{equation*}
\bgamma_\bx'(t)=\D f\big(\bgamma_\bx(t)\big),\  \bgamma_\bx(0)=\bx,
\end{equation*}
the unstable manifold for the negative gradient flow becomes the stable manifold for the gradient
flow and viceversa. Therefore, we could equivalently define the cluster associated to a mode
$\bx_0$ of the density as its stable manifold with respect to the gradient flow, that is,
$W^s_+(\bx_0)=\{\bx\in M\colon \lim_{t\to\infty}\bgamma_\bx(t)=\bx_0\big\}$. This is a precise
formulation of the notion of domain of attraction of the mode $\bx_0$, since $W^s_+(\bx_0)$
represents the set of all the points that climb to $\bx_0$ when they follow the steepest ascent
path defined by the gradient direction.

\section{Comparing population clusterings}\label{sec:cons}

Once the population goal of the clustering problem has been clearly identified as the partition of
$\mathbb R^d$ given in terms of the stable manifolds of the density modes with respect to the
gradient flow, it is natural to wonder how this population clustering can be estimated from the
data, how the error between a data-based clustering and the population clustering can be measured
and to set up a notion of consistency. Those are the goals for this section.

At this point it is worth distinguishing between two different, although closely related, concepts.
When the density $f$ is unknown, and a sample $\bX_1,\dots,\bX_n$ drawn from $f$ is given, any
procedure to obtain a data-based partition of $\mathbb R^d$ will be called a {\it data-based
clustering}. Notice, however, that when data are available the goal of most clustering procedures
is just to partition the data set, and not the whole space $\mathbb R^d$. This will be referred to
henceforth as a {\it clustering of the data}. Naturally, the ideal population clustering or any
data-based clustering immediately result in a clustering of the data, by assigning the same group
to data points belonging to the same component of the given partition of $\mathbb R^d$.

A good deal of measures to evaluate the distance between two clusterings of the data have been
proposed in the literature. The work of \cite{M07} provides both a comprehensive survey of the most
used existing measures as well as a deep technical study of their main properties, and for instance
\cite{AB73} or \cite{D81} provide further alternatives. But apart from some exceptions, as the
Hamming distance considered in \cite{BLP06}, it should be stressed that all these proposals concern
partitions of a finite set. However, since the ideal population goal of clustering introduced here
is defined as a partition of $\mathbb R^d$, it is necessary to introduce distances between two
partitions of $\mathbb R^d$ to evaluate the differences between two clusterings.

\subsection{A distance in measure between clusterings}

Let $\mathscr C$ and $\mathscr D$ be two clusterings of a probability distribution $P$, and assume
for the moment that both have the same number of clusters, say $\mathscr C=\{C_1,\dots,C_r\}$ and
$\mathscr D=\{D_1,\dots,D_r\}$. The first step to introduce a distance between $\mathscr C$ and
$\mathscr D$ is to consider a distance between sets. Surely the two distances between sets most
used in practice are the Hausdorff distance and the distance in measure; see \cite{CF10}. The
Hausdorff distance is specially useful when dealing with compact sets (it defines a metric in the
space of all compact sets of a metric space), it tries to capture the notion of physical proximity
between two sets \citep{RC03}. In contrast, the distance in measure between two sets $C$ and $D$
refers to the content of their symmetric difference $C\triangle D=(C\cap D^c)\cup(C^c\cap D)$, and
it defines a metric on the set of all measurable subsets of a measure space, once two sets
differing in a null-measure set are identified to be the same.

Although we will return to the Hausdorff distance later, our first approach to the notion
of distance between $\mathscr C$ and $\mathscr D$ relies primarily on the concept of distance in measure, and the measure involved is precisely the probability measure $P$. From a practical point of view, it does not seem so important that the clusters of a data-based partition get physically
close to those of the ideal clustering. Instead, it is desirable that the points that are incorrectly assigned do not represent a very significant portion of the distribution. Thus, it seems natural to affirm that two clusters $C\in\mathscr C$ and $D\in\mathscr D$ (resulting from different clusterings) are close
when $P(C\triangle D)$ is low.

Therefore, for two clusterings $\mathscr C$ and $\mathscr D$ with the same number of clusters, a sensible notion of distance is obtained by adding
up the contributions of the pairwise distances between their components, once they have been
relabelled so that every cluster in $\mathscr C$ is compared with its most similar counterpart in
$\mathscr D$. In mathematical terms, the distance between $\mathscr C$ and $\mathscr D$ can be
measured by
\begin{equation}\label{d1dist}
d_1(\mathscr C,\mathscr D)=\min_{\sigma\in \mathcal P_r}\sum_{i=1}^rP(C_i\triangle D_{\sigma(i)}),
\end{equation}
where $\mathcal P_r$ denotes the set of permutations of $\{1,2,\ldots,r\}$.

It can be shown that $d_1$
defines a metric in the space of all the partitions with the same number of components, once two
such partitions are identified to be the same if they differ only in a relabelling of their
components. Moreover, the minimization problem in (\ref{d1dist}) is usually known as the {\it linear sum assignment problem} in the literature of Combinatorial Optimization, and it represents a particular case of the well-known Monge-Kantorovich transportation problem. A comprehensive treatment of assignment problems can be found in \cite{BDM09}.

If a partition is understood as a vector in the product space of measurable sets, with the
components as its coordinates, then $d_1$ resembles the $L_1$ product distance, only adapted to
take into account the possibility of relabelling the components. This seems a logical choice given
the additive nature of measures, as it adds up the contribution of each distance between the partition components as
described before. However, it would be equally possible to consider any other $L_p$ distance,
leading to define
$$d_p(\mathscr C,\mathscr D)=\min_{\sigma\in \mathcal P_r}\Big\{\sum_{i=1}^rP(C_i\triangle
D_{\sigma(i)})^p\Big\}^{1/p}$$ for $p\geq1$ and also $d_\infty(\mathscr C,\mathscr
D)=\min_{\sigma\in \mathcal P_r}\max\{P(C_i\triangle D_{\sigma(i)})\colon i=1,\dots,r\}$. The
minimization problem defining $d_\infty$ is also well known, under the name of the linear
bottleneck assignement problem, and its objective function is usually employed if the interest is
to minimize the latest completion time in parallel computing \citep[see][Section 6.2]{BDM09}. In
the context of clustering, surely the $d_1$ distance seems the most natural choice among all the
$d_p$ possibilities, due to its clear interpretation.

Nevertheless, the definition of the $d_1$ distance involves some kind of redundancy, due to
the fact that $\mathscr C$ and $\mathscr D$ are partitions of $\mathbb R^d$, because the two
disjoint sets that form every symmetric difference in fact appear twice in each of the sums in
(\ref{d1dist}). More precisely, taking into account that $P(C\triangle D)=P(C)+P(D)-2P(C\cap D)$, it follows that for every $\sigma\in\mathcal P_r$
\begin{equation}\label{fact2}
\sum_{i=1}^rP(C_i\triangle D_{\sigma(i)})=2-2\sum_{i=1}^r P(C_i\cap D_{\sigma(i)}).
\end{equation}
To avoid this redundancy, the eventual suggestion to measure the distance between $\mathscr C$ and
$\mathscr D$, based on the set distance in $P$-measure, is $d_P(\mathscr C,\mathscr D)=\frac12 d_1(\mathscr C,\mathscr D)$.

If the partitions $\mathscr C$ and $\mathscr D$ do not have the same number of clusters, then as many empty set components as needed are added so that both partitions include the same number of components, as in \cite{CDGH06}, and the distance between the extended partitions is computed as before. Explicitly, if $\mathscr C=\{C_1,\dots,C_r\}$ and $\mathscr D=\{D_1,\dots,D_s\}$ with $r<s$ then, writing $C_i=\varnothing$ for $i=r+1,\dots,s$, we set
$$d_P(\mathscr C,\mathscr D)=\frac12\min_{\sigma\in\mathcal P_s}\sum_{i=1}^sP(C_i\triangle D_{\sigma(i)})=\frac12\min_{\sigma\in\mathcal P_s}\Big\{\sum_{i=1}^rP(C_i\triangle D_{\sigma(i)})+\sum_{i=r+1}^sP(D_{\sigma(i)})\Big\}.$$
Thus, the term $\sum_{i=r+1}^sP(D_{\sigma(i)})$ can be interpreted as a penalization for unmatched probability mass.

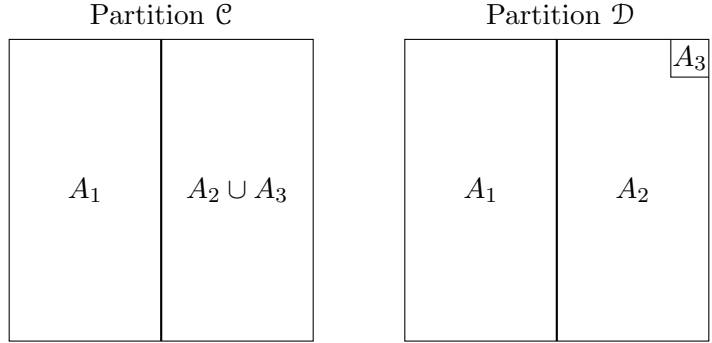
\begin{figure}[t]\centering
\begin{tabular}{cc}
Partition $\mathscr C$&\qquad Partition $\mathscr D$\\
\begin{tikzpicture}[scale=0.5]
    \draw (-4,-4) -- (4,-4) -- (4,4) -- (-4,4) -- cycle;
    %\fill [fill=green!60] (0,-4) -- (4,-4) -- (4,4) -- (0,4) -- cycle;
    \draw [thick] (0,-4) -- (0,4);
    %\draw (0,0) node[] {$\bullet$};
    %\draw (0,0) node[below] {\small $\bx_0$};
    %\draw (-1.5,0) node[] {$\bullet$};
    \draw (-2,0) node {$A_1$};
    %\draw (1.5,0) node[] {$\bullet$};
    \draw (2,0) node {$A_2\cup A_3$};
\end{tikzpicture}&\qquad\begin{tikzpicture}[scale=0.5]
    \draw (-4,-4) -- (4,-4) -- (4,4) -- (-4,4) -- cycle;
    %\fill [fill=green!60] (0,-4) -- (4,-4) -- (4,4) -- (0,4) -- cycle;
    \draw [thick] (0,-4) -- (0,4);
    %\draw (0,0) node[] {$\bullet$};
    %\draw (0,0) node[below] {\small $\bx_0$};
    %\draw (-1.5,0) node[] {$\bullet$};
    \draw (-2,0) node {$A_1$};
    %\draw (1.5,0) node[] {$\bullet$};
    \draw (2,0) node {$A_2$};
    \draw (3,4) -- (3,3) -- (4,3);
    \draw (3.5,3.5) node {$A_3$};
\end{tikzpicture}
\end{tabular}
\caption{Two partitions of the unit square that do not differ much if $A_3$ has low probability.} \label{fig:5}
\end{figure}

The idea is that two partitions as those shown in Figure \ref{fig:5} do not differ much if $A_3$
has low probability, even if they do not have the same number of clusters. For the partitions in
Figure \ref{fig:5}, denote $\mathscr C=\{C_1,C_2\}$ and $\mathscr D=\{D_1,D_2,D_3\}$ with
$C_1=D_1=A_1$, $C_2=A_2\cup A_3$, $D_2=A_2$, $D_3=A_3$, and assume that $P(A_1)=0.5$, $P(A_2)=0.45$
and $P(A_3)=0.05$. Then, it can be shown that $d_P(\mathscr C, \mathscr D)=0.05$. In (\ref{d1dist})
every cluster of $\mathscr C$ is matched to some cluster in $\mathscr D$, depending on the
permutation for which the minimum is achieved. When $\mathscr C$ has less clusters than $\mathscr
D$, some of the components of $\mathscr D$ will be matched with the empty set, indicating that they
do not have an obvious match in $\mathscr C$, or that they are unimportant. In the previous
example, the minimum is achieved when $C_1$ is matched with $D_1$, $C_2$ with $D_2$ and $D_3$ is
matched with the empty set.

It is interesting to note that $d_P(\mathscr C,\mathscr D)$ can be estimated in a natural way by
replacing $P$ with the empirical measure based on the data $\bX_1,\dots,\bX_n$, leading to
$$\widehat{d_P}(\mathscr C,\mathscr D)=\frac1{2n}\min_{\sigma\in\mathcal P_s}\Big\{\sum_{i=1}^r\sum_{j=1}^nI_{C_i\triangle D_{\sigma(i)}}(\bX_j)+\sum_{i=r+1}^s\sum_{j=1}^nI_{D_{\sigma(i)}}(\bX_j)\Big\}$$
where $I_A$ denotes the indicator function of the set $A$. When $r=s$, it follows from (\ref{fact2}) that an alternative expression for $d_P(\mathscr C,\mathscr D)$ is
$$d_P(\mathscr C,\mathscr D)=1-\max_{\sigma\in\mathcal P_r}\sum_{i=1}^rP(C_i\cap D_{\sigma(i)})$$
and, therefore, its sample analogue
$$\widehat{d_P}(\mathscr C,\mathscr D)=1-\frac1n\max_{\sigma\in\mathcal P_r}\sum_{i=1}^r\sum_{j=1}^nI_{C_i\cap D_{\sigma(i)}}(\bX_j),$$
coincides with the so-called classification distance between two clusterings of the data, whose
properties are explored in \citet{M05,M07,M12}. For $r<s$, however, $\widehat{d_P}$ differs from
the classification distance (which does not include the penalty term), but it corresponds exactly
with the transfer distance, studied in detail in \cite{CDGH06} \citep[see also][]{D08}. Extending
the properties of the transfer distance to its population counterpart suggests an interpretation of
$d_P(\mathscr C,\mathscr D)$ as the minimal probability mass that needs to be moved to transform
the partition $\mathscr C$ into $\mathscr D$, hence the connection with the optimal transportation
problem.

\subsection{A distance between sets of sets}

An alternative notion of distance between two clusterings based on the Hausdorff metric has been kindly suggested by Professor Antonio Cuevas, noting that precisely this distance was used in \cite{P81} to measure the discrepancy between the set of sample $K$-means and the set of population $K$-means. If $(X,\rho)$ is a metric space and $A,B\subseteq X$ are two non-empty subsets of $X$, the Hausdorff distance between $A$ and $B$ is defined as
$$d_H(A,B)=\max\Big\{\sup_{a\in A}\inf_{b\in B}\rho(a,b),\,\sup_{b\in B}\inf_{a\in A}\rho(a,b)\Big\}$$
or, equivalently, as
$$d_H(A,B)=\inf\{\varepsilon>0\colon A\subseteq B^{\varepsilon}\text{ and }B\subseteq A^{\varepsilon}\Big\},$$
where $A^\varepsilon=\bigcup_{a\in A}\{x\in X\colon \rho(x,a)\leq\varepsilon\}$, and $B^\varepsilon$ is defined analogously.

In the context of clustering, take $X$ to be the metric space consisting of all the sets of $\mathbb R^d$ equipped with the distance $\rho(C,D)=P(C\triangle D)$, once two sets with $P$-null symmetric difference have been identified to be the same. Then any two clusterings $\mathscr C=\{C_1,\dots,C_r\}$ and $\mathscr D=\{D_1,\dots,D_s\}$ can be viewed as (finite) subsets of $X$, and therefore the Hausdorff distance between $\mathscr C$ and $\mathscr D$ is defined as
\begin{align*}
d_H(\mathscr C,\mathscr D)&=\max\Big\{\max_{i=1,\dots,r}\min_{j=1,\dots,s}P(C_i\triangle D_j),\,\max_{j=1,\dots,s}\min_{i=1,\dots,r}P(C_i\triangle D_j)\Big\}\\
&=\inf\{\varepsilon>0\colon \mathscr C\subseteq \mathscr D^{\varepsilon}\text{ and }\mathscr D\subseteq \mathscr C^{\varepsilon}\Big\}.
\end{align*}
To express it in words, $d_H(\mathscr C,\mathscr D)\leq\varepsilon$ whenever for every $C_i\in\mathscr C$ there is some $D_j\in\mathscr D$ such that $P(C_i\triangle D_j)\leq\varepsilon$ and viceversa.

The Hausdorff distance can be regarded as a uniform distance between sets. It is not hard to show, using standard techniques from the Theory of Normed Spaces, that when $r=s$ we have
$$d_H(\mathscr C,\mathscr D)\leq 2\, d_P(\mathscr C,\mathscr D)\leq r\, d_H(\mathscr C,\mathscr D).$$
However, when $r<s$ the distance $d_H$ can be more demanding than $d_P$, meaning that both
partitions have to be really close so that their Hausdorff distance results in a small value. For
instance, it can be checked that for the two clusterings of the previous example, shown in Figure
\ref{fig:5}, the Hausdorff distance between them is $d_H(\mathscr C,\mathscr D)=0.45$, mainly due
to the fact that $C_2$ and $D_3$ are far from each other, since $P(C_2\triangle D_3)=P(A_2)=0.45$.
A clear picture of the difference between $d_H$ and $d_P$ is obtained by arranging all the
component-wise distances $P(C_i\triangle D_j)$ into an $r\times s$ matrix. Then, the Hausdorff
distance is obtained by computing all the row-wise and column-wise minima and taking the maximum of
all of them. In contrast, for the distance in $P$-measure the first step when $r<s$ is to add $s-r$
row copies of the vector $(P(D_1),\dots,P(D_s))$ to the matrix of component-wise distances, and
then compute the distance in $P$-measure as half the minimum possible sum obtained by adding up a
different element in each row. As a further difference, note that the Hausdorff distance does not
involve a matching problem; instead, this distance is solely determined by the two components that
are furthest from each other.

Obviously, a sample analogue is also obtained in this case by replacing $P$ for the empirical
probability measure, leading to
\begin{align*}
\widehat d_H(\mathscr C,\mathscr D)=\frac1n\max\Big\{\max_{i=1,\dots,r}\min_{j=1,\dots,s}\sum_{k=1}^nI_{C_i\triangle D_j}(\bX_k),\,
\max_{j=1,\dots,s}\min_{i=1,\dots,r}\sum_{k=1}^nI_{C_i\triangle D_j}(\bX_k)\Big\},
\end{align*}
which seems not to have been considered previously as a distance between two clusterings of the
data.

\subsection{Consistency of data-based clusterings}

As indicated above, a data-based clustering is understood as any procedure that induces a partition
$\widehat{\mathscr C}_n$ of $\mathbb R^d$ based on the information obtained from a sample
$\bX_1,\dots,\bX_n$ from a probability distribution $P$. Once we have a precise definition for the
ideal population clustering goal $\mathscr C$, a data-based clustering $\widehat{\mathscr C}_n$ is
said to be consistent if it gets closer to $\mathscr C$ as the sample size increases. Formally, if
$d(\widehat{\mathscr C}_n,\mathscr C)\to0$ as $n\to\infty$ for some of the modes of stochastic
convergence (in probability, almost surely, etc), where $d$ represents one of the distances between
clusterings defined above, or any other sensible alternative.

A plug-in strategy to obtain data-based clusterings naturally arises when the density $f$ is
replaced with a smooth estimator $\hat f$, and so the data-based clusters are defined as the
domains of attraction of the modes of $\hat f$. Obvious candidates for the role of $\hat f$ include
kernel density estimators for the nonparametric setup or mixture model density estimators
in a parametric context.

This is a very simple approach, that involves to some extent estimating the density function to
solve the clustering problem (unsupervised learning), and according to \citet[][p. 21]{vL04} it may
not be a good idea because density estimation is a very difficult problem, especially in high
dimensions. A similar situation is found in the study of classification (supervised learning),
where the optimal classifier, the Bayes rule, depends on the regression function of the random
labels over the covariates. Here, even if classification can be proved to be a problem easier than
regression, nevertheless regression-based algorithms for classification play an important role in
the development of supervised learning theory \cite[see][Chapter 6]{DGL96}.

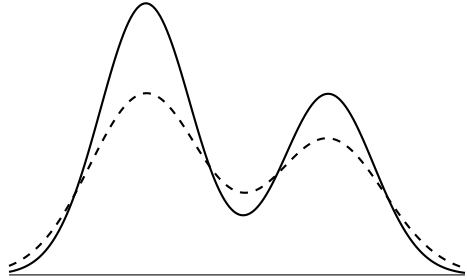
\begin{figure}[t]\centering
\begin{tikzpicture}[scale=0.6]
    \def\normalmixt{\x,{6*exp(-((\x)^2)/2)+4*exp(-((\x-4)^2)/2)}}
    \def\normalmix{\x,{4*exp(-((\x)^2)/3)+3*exp(-((\x-4)^2)/3)}}
    \draw[color=black,domain=-3:7,samples=200,thick] plot (\normalmixt) node[right] {};
    \draw[color=black,domain=-3:7,samples=200,dashed,thick] plot (\normalmix) node[right] {};
    \draw (-3,0) -- (7,0) node[right] {};
\end{tikzpicture}
\caption{Two density functions that are not close but induce exactly the same clustering.}
\label{fig:6}
\end{figure}

Along the same lines, Figure \ref{fig:6} illustrates why we should not completely discard density estimation as an intermediate step for clustering. Figure \ref{fig:6} shows a typical situation where the solid line is the true density and the dashed line is a kernel density estimator, since an expansion of its pointwise bias shows that on average the kernel estimator underestimates the maxima and overestimates the minima \cite[][p. 21]{WJ95}. But even if the two density functions are not really close in any global sense, they produce exactly the same clusterings of $\mathbb R$, so this also seems to suggest that clustering should be easier than density estimation. Studying consistency and error rates for data-based clusterings obtained from density estimators surely represents an interesting open problem.

\section{Practical algorithms}

Apart from the theoretical considerations in the previous sections, for practical purposes it is necessary to provide algorithms that help to compute population or data-based clusterings. As noted in Section \ref{sec:1}, even if the density $f$ is fully known this may become an onerous task for some methodologies, as
for instance those based on level sets, for which smart implementations are needed to compute the connected components of the level sets or the cluster tree \citep{CFF01,SN10}.

With the concept of cluster as the domain of attraction of the density modes, the natural approach
to compute the partition induced by a density is to use the so-called mode seeking algorithms. Surely the most
popular algorithm in this class is the mean shift algorithm, initially introduced by \cite{FH75}
and subsequently recast by \cite{Ch95} and \cite{CM02} to highlight several applications in
Engineering. The mean shift algorithm is a iterative procedure which, at every step, moves the
point obtained in the previous iteration to a location of higher density, thus producing a
convergent sequence that transports any initial value to a local maximum of the density.
Explicitly, any initial point $\by_0$ is transformed recursively to obtain a sequence defined by
\begin{equation}\label{eq:meanshift}
\by_{j+1}=\by_j+\mat A {\D f}(\by_j)\big/f(\by_j),
\end{equation}
where $\mat A$ is a positive-definite matrix chosen to ensure convergence. Notice that if $\mat A$
is a positive multiple of the identity matrix then the sequence approximately follows the steepest
ascent path, since at every step the shift is made along the gradient direction. In this case, the
mean shift algorithm is easily recognized as a variant of the classical gradient ascent algorithm
used for numerical maximization, but employing the normalized gradient (the gradient divided by
$f$) to accelerate convergence in low-density zones.

Since this algorithm can be applied to any initial point, it effectively produces a partition of
$\mathbb R^d$. However, when a sample $\bX_1,\dots,\bX_n$ is available it is common to use every
sample point as initial value to obtain a clustering of the data. On the other hand, if the sample is employed to construct a smooth estimator $\hat f$ of the
density, then Equation (\ref{eq:meanshift}) with $f$ replaced by $\hat f$ allows to produce a
data-based clustering (hence, also a clustering of the data as noted above).

Moreover, for certain types of density estimators there exist closed forms for the matrix $\mat A$ that ensure the
convergence of the algorithm to a local maximum. In general, suppose the density $f$ corresponds to
a normal mixture model, so that
\begin{equation}\label{nmd}
f(\bx)=\sum_{\ell=1}^L\pi_\ell\phi(\bx|\bmu_\ell,\bSigma_\ell),
\end{equation}
where
$\phi(\bx|\,\bmu,\bSigma)=|2\pi\bSigma|^{-1/2}\exp\{-\frac12(\bx-\bmu)^\top\bSigma^{-1}(\bx-\bmu)\}$
refers to the density of the normal distribution with mean vector $\bmu$ and covariance matrix
$\bSigma$, and $\pi_1,\dots,\pi_L$ are nonnegative quantities that sum to one. Then, \cite{CP07}
showed that the iterative algorithm defined by (\ref{eq:meanshift}) is convergent for every initial
value $\by_0$ by taking $\mat A\equiv\mat A_j=\bSigma_*(\by_j)$, where $\bSigma_*(\by)=\{\sum_{\ell=1}^L\alpha_\ell(\by)\bSigma_\ell^{-1}\}^{-1}$ is a weighted harmonic mean of the covariance matrices of the components, with weights
$\alpha_\ell(\by)=\pi_\ell\phi(\by|\bmu_\ell,\bSigma_\ell)/f(\by)$ given by the posterior odds of the $\ell$th component of the mixture after observing $\by$.

\begin{figure}[t]\centering
\includegraphics[width=0.45\textwidth]{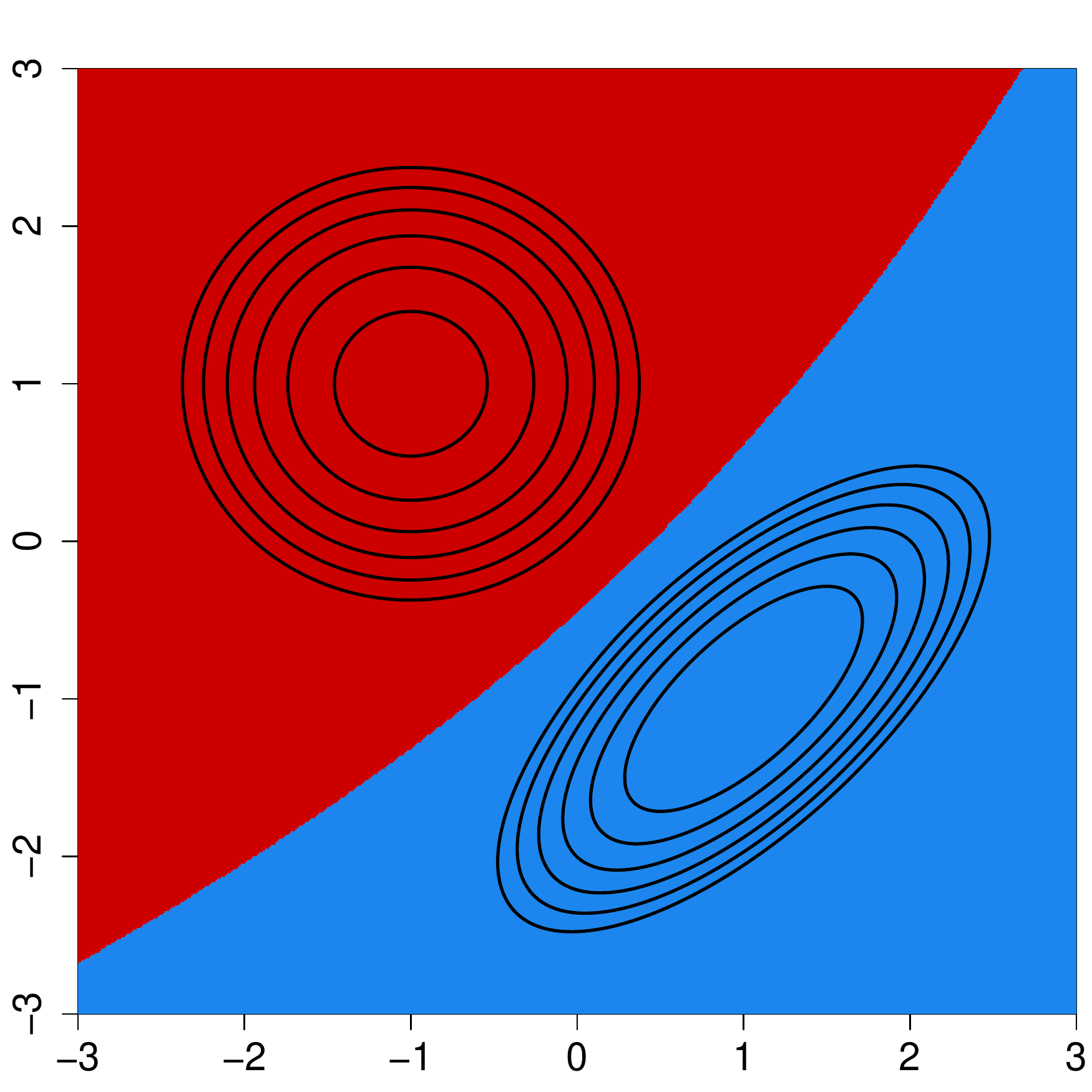}\hspace{0.05\textwidth}\includegraphics[width=0.45\textwidth]{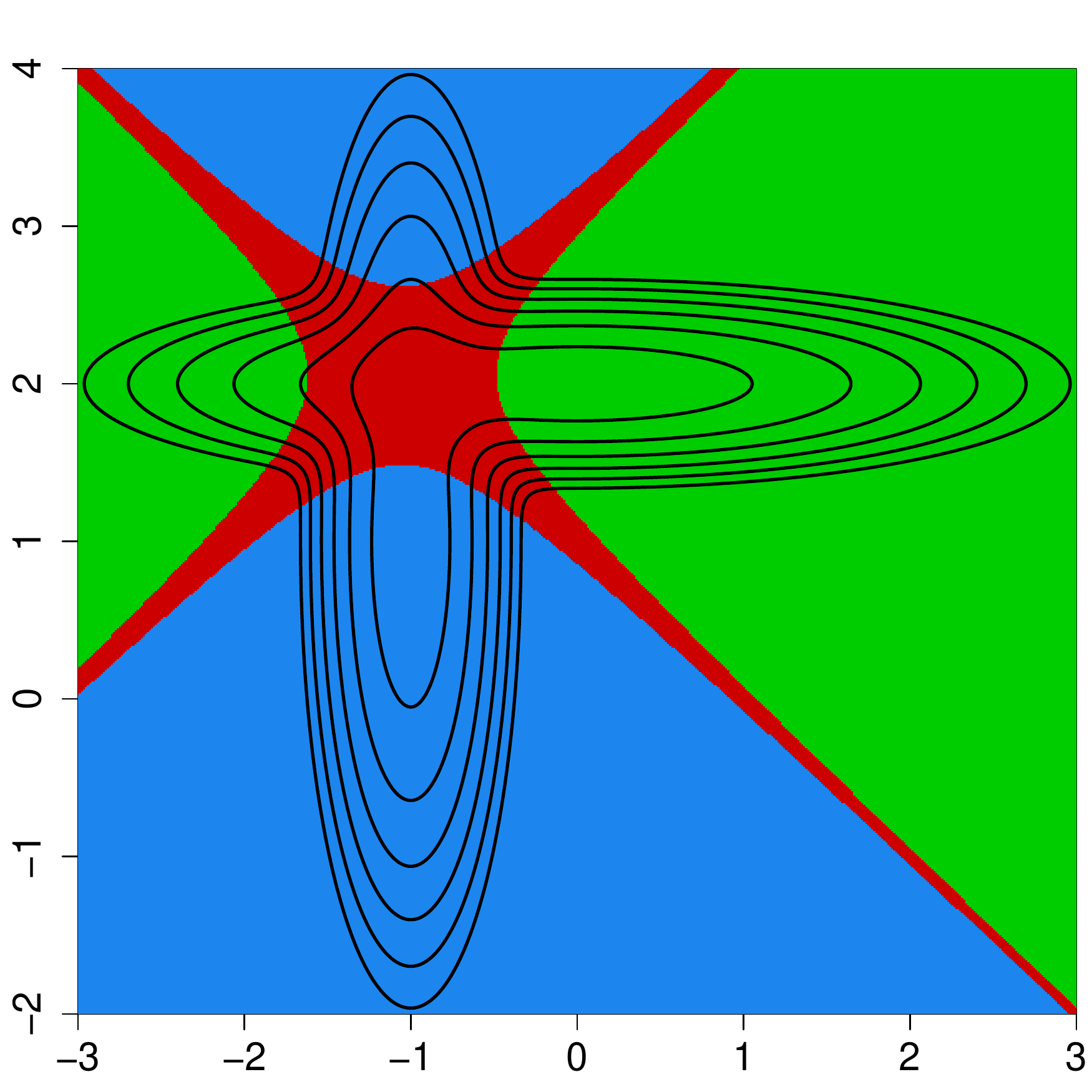}\\
\includegraphics[width=0.45\textwidth]{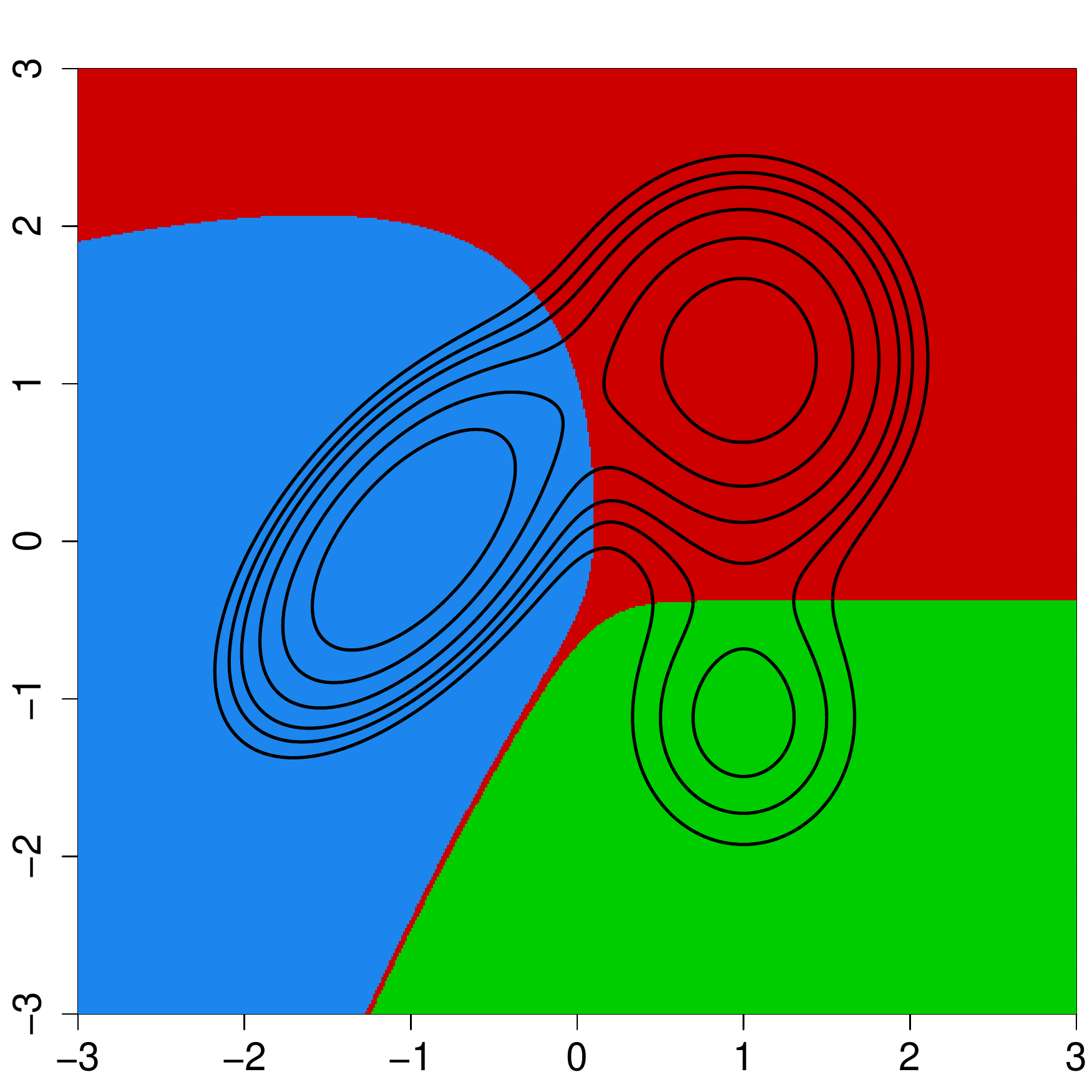}\hspace{0.05\textwidth}\includegraphics[width=0.45\textwidth]{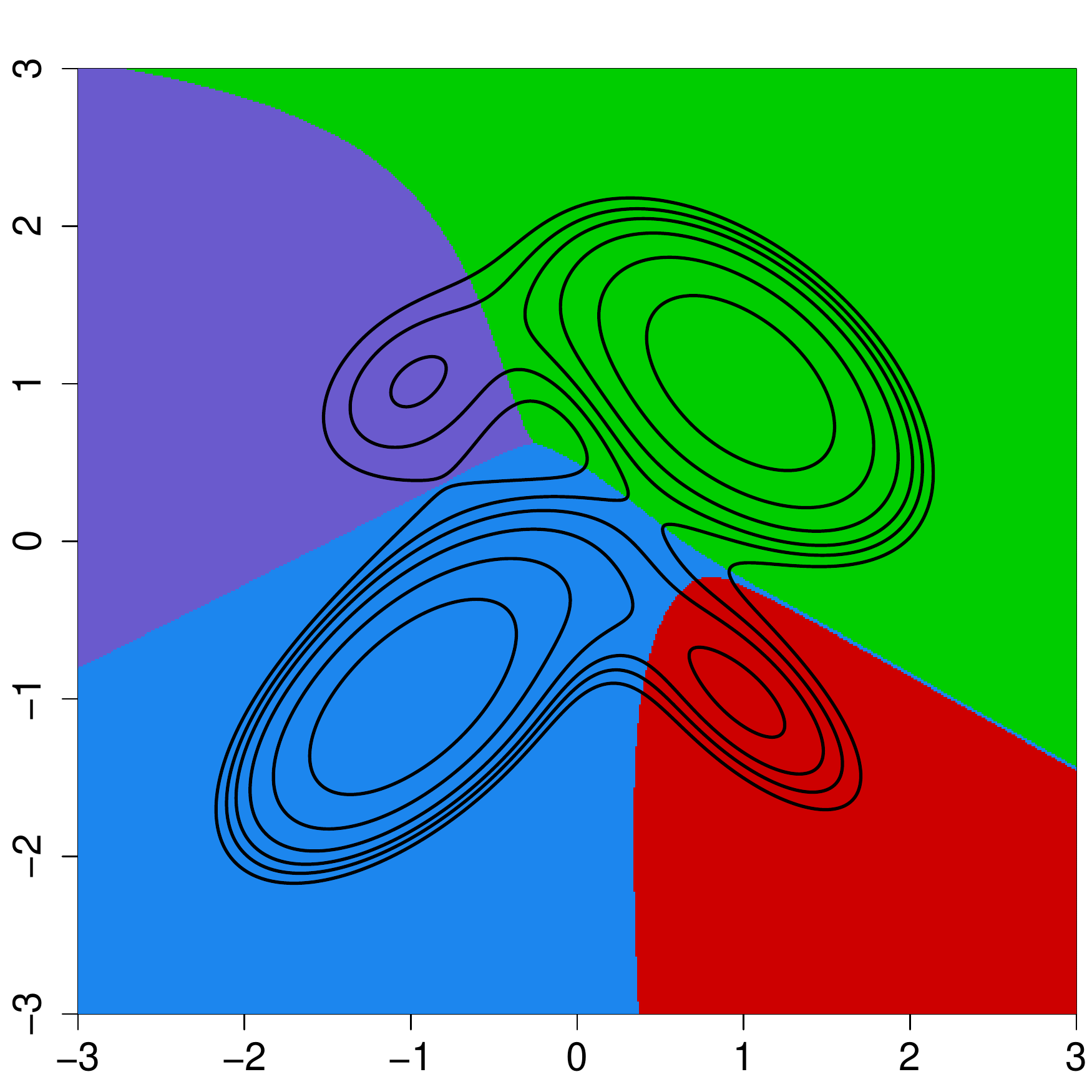}
\caption{Examples showing the ideal population clustering for several bivariate normal mixture densities, calculated using the mean shift algorithm.}
\label{fig:7}
\end{figure}

Hence, in the normal mixture model where the density estimator $\hat f$ is obtained by replacing the parameters in (\ref{nmd}) with estimates $\hat\pi_\ell,\hat\bmu_\ell,\hat\bSigma_\ell$, the mean shift algorithm with $\mat A=\hat{\bSigma}_*(\by_j)$ (i.e., the previous definition, but also with the parameters replaced by estimates) provides a computationally simple way to obtain the corresponding data-based clustering. Some examples of the partitions thus obtained, for normal mixture densities taken from \cite{WJ93} and \cite{RL05}, are shown in Figure \ref{fig:7}.

Furthermore, suppose that a nonparametric kernel density estimator $\hat f_\bH$ is employed in (\ref{eq:meanshift}), namely $\hat f_{\bH}(\bx)=n^{-1}\sum_{i=1}^nK_\bH(\bx-\bX_i)$, where the kernel $K$ is a fixed density function, $\bH$ is a positive-definite bandwidth matrix and $K_\bH(\bx)=|\bH|^{-1/2}K(\bH^{-1/2}\bx)$. If the kernel $K$ is taken as the normal density, then $\hat f_\bH(\bx)=\sum_{i=1}^nn^{-1}\phi(\bx|\bX_i,\bH)$, and so its structure as a normal mixture density allows to conclude as before that the mean shift algorithm is convergent for any initial point by just taking $\mat A=\bH$. More generally, proceeding as in the proof of Theorem 1 in \cite{CM02} it can be showed that convergence is guaranteed with $\mat A=\bH$ for a general kernel of the form $K(\bx)=k(\bx^\top\bx)$, as long as the profile $k:[0,\infty)\to\mathbb R$ is convex and monotonically decreasing.

Many variants of the mean shift algorithm have been explored in the literature, expanding in different directions. In the nonparametric context, the iterative procedure (\ref{eq:meanshift}) with $\mat A=\bH$ as above can be explicitly written as $\by_{j+1}=\sum_{i=1}^nW_i(\by_j)\bX_i$, hence taking the form of a weighted mean of the observations, with the weight functions defined as
\begin{equation*}\label{eq:wms}
W_i(\by)=\frac{G_\bH(\by-\bX_i)}{\sum_{i=1}^nG_\bH(\by-\bX_i)},
\end{equation*}
where the kernel $G(\bx)\propto g(\bx^\top\bx)$ is based on the positive profile $g=-k'$. \cite{SKK07} provided an alternative formulation of (\ref{eq:meanshift}) by noting that
\begin{equation}\label{eq:ms2}
\by_{j+1}=\mathop{\rm argmin}_{\by\in\mathbb R^d}\sum_{i=1}^n\|\bX_i-\by\|_2^2\,G_\bH(\bX_i-\by_j),
\end{equation}
where $\|\bx\|_2=(\bx^\top\bx)^{1/2}$ denotes the usual Euclidean distance, and proposed the medioid-shift algorithm by restricting the minimization in (\ref{eq:ms2}) to the set of sample points, thus reducing the computation time since only one iteration of the algorithm has to be computed per data point because $\by_{j+1}=\bX_k$ for some $k\in\{1,\dots,n\}$. Accelerating the mean shift algorithm is in fact one of the expansion directions that has generated much research. A (clearly non-exhaustive) list of references on this topic would include \cite{CP06}, \cite{Wal07} and \cite{VS08}.

Other variants of the mean shift algorithm are oriented towards robustness, by considering median shifts instead of mean shifts. \cite{WQZ07} proposed the algorithm {\sc clues} that at each iteration shifts the previous point to the coordinate-wise (i.e., marginal) median of its $K$-nearest neighbors and \cite{SAS09} explored other possibilities using the different concepts of multivariate median related to depth measures. A further refinement of {\sc clues}, with improved computational efficiency, and also based on local medians is the {\sc attractors} algorithm of \cite{PVZ12}. In fact, it would be possible to modify (\ref{eq:ms2}) to obtain shift algorithms based on other norm-based local medians, which would be the local analogues to the global proposals studied in \cite{DR99} \citep[see also][]{Sm90}.

Indeed, a geometric median shift algorithm over Riemannian manifolds was introduced in \cite{WH10}, and its applications to clustering were explored in \cite{WHW13}. In this respect, the adaptation of shift algorithms to manifolds does represent a third major avenue of active research on the topic. For instance, among some other recent papers dealing with this subject we could cite \cite{SM09}, \cite{OE11} and \cite{GPVW12}.

\section{Conclusions and directions for further research}

At the time of comparing different clustering procedures it is necessary to have a ``ground truth",
or population goal, that represents the ideal clustering to which the clustering algorithms should
try to get close. Sometimes this ideal population clustering is not so easy to specify and it
depends on the notion of cluster in which the researcher is interested. Here, the population goal
for modal clustering is accurately identified, making use of some tools from Morse theory, as the
partition of the space induced by the domains of attraction of the local maxima of the density
function.

This definition needs the probability density to be smooth to a certain degree, specifically it
must be a Morse function. It would be appealing to extend this notion to density functions that are
not Morse functions, meaning either that they are smooth but have degenerate critical points or
even that they are not differentiable. To treat the first case, it might be useful to resort to the
theory of singularities of differential mappings, which is exhaustively covered in the book by
\cite{AGLV98}, for instance. On the other hand, the study of the non-smooth case might start from
\cite{APS97}, where Morse theory for piecewise smooth functions is presented. Here, the key role
would be played by the subgradient, which generalizes the concept of gradient for non-smooth
functions.

Another interesting challenge consists of studying the choice of the parameters for the density estimators (the bandwidth for kernel estimators, the mixture parameters for mixture model estimators) that minimize the distance between the corresponding data-based clustering and the true population clustering, as measured by any of the distance between clusterings discussed in Section \ref{sec:cons}. Or, perhaps even better, to develop methods aimed to perform modal clustering that do not necessarily rely on a pilot density estimate, as it happens for classification methods, which are not necessarily based on a regression estimate.

On the algorithmic side, the alternative formulation (\ref{eq:ms2}) of the mean shift algorithm by \cite{SKK07} surely deserves further attention. Its resemblance to locally weighted regression methodology \citep{RW94}, identifying somehow mean shift as a locally constant adjustment, suggests that perhaps a local linear mean shift could improve the performance of clustering based on kernel density estimation, in the same way as local linear regression improves over the Nadaraya-Watson (locally constant) estimator.

\bigskip

\noindent{\bf Acknowledgments.} This work has been partially supported by Spanish Ministerio de
Ciencia y Tecnolog\'{\i}a project MTM2010-16660. The author wishes to thank Prof. Antonio Cuevas
from Universidad Aut\'onoma de Madrid as well as Prof. Ricardo Faro from Universidad de Extremadura
for insightful conversations and suggestions concerning the material of Section \ref{sec:cons}. The paper by \cite{RL05}, in which interesting connections between Morse theory and the topography of multivariate normal mixtures are illustrated, was thought-provoking enough to inspire part of this paper.

\bibliographystyle{apalike}
%\bibliography{../../../../bibliography/biblio} \end{document}

\end{document}